\documentclass[11pt,a4paper]{article}
\usepackage{amsfonts,amsmath,amstext,amssymb}
\usepackage{theorem}    %para cambiar el estilo de los teoremas
\usepackage{verbatim}

\def\R{\mathbb R}     % numeros reales
\newtheorem{lema}{Lemma}

\newtheorem{teor}[lema]{\bf Theorem}
\newtheorem{coro}[lema]{\bf Corollary}

\newtheorem{rema}[lema]{\bf Remark}

\theorembodyfont{\rmfamily}

\title{Rigidity and non-existence of spacelike submanifolds with causal mean curvature vector field  in spacetimes and the Cauchy problem in General Relativity}

\author{Juan A. Aledo${}^{a}$, Rafael M. Rubio${}^{b}$ and Juan J. Salamanca${}^{c}$ \\[6mm]
${}^a$ Departamento de Matem\'aticas, E.S.I. Inform\'atica, \\[0.5mm] Universidad de
Castilla-La Mancha, 02071 Albacete, Spain,\\ E-mail\textup{:
\texttt{juanangel.aledo@uclm.es}} \\[3mm]
${}^b$ Departamento de Matem\'aticas, Facultad de Ciencias, \\[0.5mm] Universidad de
C\'ordoba, 14071 C\'ordoba, Spain,\\[0.5mm] E-mails\textup{: \texttt{rmrubio@uco.es}, \texttt{jjsalamanca@uco.es}}\\[3mm]
${}^c$ Departamento de Estad\'{\i}stica, E.S.I. Inform\'atica, \\[0.5mm] Universidad de
Oviedo, 33003 Gij\'on, Spain,\\ E-mail\textup{:
\texttt{salamancajuan@uniovi.es}} \\[3mm]}

\date{}

\begin{document}

\maketitle

\thispagestyle{empty}

\begin{abstract}
New general results of non-existence and rigidity of spacelike
submanifolds immersed in a spacetime, whose mean curvature is a time-oriented causal vector
field,  are given. These results hold for a wide class of
spacetimes which includes globally hyperbolic, stationary,
conformally stationary and pp-wave spacetimes, among others.
Moreover, applications to the Cauchy problem in General
Relativity, are presented. Finally, in the case of hypersurfaces, we also obtain significant
consequences in Geometrical Analysis, solving new Calabi-Bernstein
and Dirichlet problems on a Riemannian manifold. \vspace*{2mm}
\end{abstract}


\vspace*{-2mm}

\section{Introduction}
In addition to geometric interest of the study of spacelike submanifolds, with causal mean curvature vector field in Lorentzian manifolds, the importance from a physical point of view is also  indisputable. So, for example,  trapped surfaces
 are nowadays a very intense
field of study.

As it is said in \cite{mars},
their applications in General Relativity are
ubiquitous: 

\begin{quote} {\it
the development of the singularity theorems, the general analysis of gravitational collapse and formation of trapped surfaces and black holes, the study of the cosmic censorship hypothesis and the related Penrose inequality], or the numerical analysis of the Cauchy development of apparently innocuous initial data}. 

\end{quote}

The concept of closed (compact without boundary) trapped surfaces,
was originally formulated by Penrose in 1965 for the case of
2-dimensional spacelike surfaces in $4$-dimensional spacetimes in
terms of the signs or the vanishing of the so-called null expansions
(see \cite{Pe}, \cite{HP}), and it has remained that way for many
years. Note that  the Penrose's definition implies that a closed
trapped surface is a topological sphere. In recent years it has
become clear that this concept is related to the causal orientation
of the mean curvature vector field of the submanifold, which
provides a better characterization of the trapped surfaces and
allows their generalization to 2-codimensional spacelike
submanifolds of arbitrary dimension $n$ (see \cite{Sen} for more
details). Recall that an embedded  2-codimensional spacelike
submanifold $S$ is called a \emph{future} (resp. \emph{past})
\emph{trapped surface} if its mean curvature vector field $\vec{H}$
is timelike and future pointing (resp. past pointing). When
$\vec{H}$ is causal and future pointing all over the spacelike
submanifold and timelike at least at a point of $S$, the submanifold
is called \emph{nearly future trapped}, and correspondingly for
\emph{nearly past trapped}. When $\vec{H}$ is causal and
future-pointing everywhere, and non-zero at least at a point of $S$,
the submanifold is said to be \emph{weakly future trapped},
similarly for \emph{weakly past trapped}. Finally, the submanifold
$S$ is called \emph{marginally future trapped} if $\vec{H}$ is
ligthlike and future pointing all over $S$ and non-zero at least at
a point of $S$, and analogously for the past case. The case
$\vec{H}\equiv \vec{0}$ corresponds to extremal or symmetric
submanifolds.

There are many recent works about trapped surfaces (this understood
in a broad sense, i.e. marginally, nearly, etc). These works deal
with different aspects of the surfaces: geometrical properties,
rigidity, representation, non-existence results, classification in
some spacetimes, causality properties, etc (see for instance
\cite{AGM}, \cite{IB}, \cite{Fl}, \cite{Be}, \cite{Be}, \cite{Mu} \cite{Al},
\cite{Cr} and references therein).

In this work  we consider  spacelike
submanifolds with arbitrary codimension (including hypersurfaces)  and weak assumptions on the
mean curvature vector field, which are immersed in several relevant families of spacetimes. Notice that the study of spacelike submanifolds of arbitrary dimension and causal mean curvature is interesting in General Relativity (see for instance \cite{Ga}).  More specifically, we study immersed compact
spacelike submanifolds with time-oriented causal mean curvature vector field
where, as  usual, we take the zero vector to be null or lightlike
(future and past), but not timelike or spacelike. This paper is mathematical in nature, our study is  intrinsic, i.e., free of coordinates, so it is given in terms of modern differential geometry.

 On the other hand,
when a smooth causal vector field is not equal to the zero vector at any
point, we will say that the vector field is \emph{strictly causal}.
A causal vector field is said to be \emph{time-oriented} if it is
future causal or past causal at every point. Note that every
strictly causal smooth vector field is time-oriented.

Hence, our study is very general and it can be applied to a wide
class of relevant spacetimes. After developing several suitable
geometrical equations in section \ref{s3}, we obtain the following
general result (see Theorem \ref{r1}).

\begin{quote} {\it  (a) Let $(M,\overline{g})$ be a spacetime, which admits a future causal
vector field $X$ such that the $2$-covariant symmetric tensor
$\mathcal{L}_X \overline{g}$ is positive semi-definite (resp.
negative semi-definite) and positive definite (resp. negative definite) at some point on spacelike vectors. Then there are no
compact spacelike submanifold in $(M,\overline{g})$ with future (resp.
past) causal mean curvature, in particular there are no extremal spacelike submanifold.

\vspace{2mm}

\noindent(b) Let $(M,\overline{g})$ be a spacetime, which admits a future strictly causal
vector field $X$ such that the $2$-covariant symmetric tensor
$\mathcal{L}_X \overline{g}$ is positive semi-definite (resp.
negative semi-definite) on spacelike vectors. Then there are no
compact spacelike submanifolds in $(M,\overline{g})$ with future (resp.
past) causal mean curvature, whose mean curvature vector field
$\vec{H}$ is not zero at least at a point.}
\end{quote}

On other hand, when the spacetime admits some infinitesimal symmetry
in form of a strictly causal Killing vector field, as a direct
consequence of the previous result we obtain:

\begin{quote} {\it In a spacetime which admits a strictly causal Killing vector field,
there exists no compact spacelike submanifold with future or past
causal mean curvature vector field $\vec{H}$ except, perhaps, with
$\vec{H}\equiv \vec{0}$.

\noindent In particular, in a stationary or pp-wave spacetime, there
exists no compact spacelike submanifold with time-oriented causal
mean curvature vector field except, perhaps, with $\vec{H}\equiv
\vec{0}$.}
\end{quote}

This last result provides an interesting application to the Cauchy
problem in General Relativity (Theorem \ref{cauchy}):

\begin{quote} {\it Let $(S,g,A)$ be an initial data set for the Cauchy problem of the Einstein's equation. Assume that the Riemannian manifold  $S$
admits a non-minimal compact submanifold $P$, whose mean curvature
vector field $\vec{h}$ in $S$ satisfies
\[
\left\|\vec{h}\right\|< \left|\mathrm{trace}_{|_P} A\right| \, .
\]
Then, a development spacetime $(M,\tilde{g})$ of this data set
cannot be  a stationary spacetime. In fact, such development cannot
be a spacetime with a strictly causal Killing vector field.}
\end{quote}

In a more general framework, when the spacetime admits a conformal
infinitesimal symmetry, other interesting results for spacelike
submanifolds are described in Section \ref{main} (see, for instance
Corollaries \ref{con1} and \ref{con2}). Notice that our results in
this last direction hold in the important class of conformally
stationary spacetimes.

\vspace{2mm}

In Section \ref{s5}, we obtain again  results relitive to the Cauchy problem, so we give several theorems of non-existence of
compact spacelike submanifolds with time-oriented causal mean curvature vector
field in development spacetimes from certain initial data set. For
example, we obtain the result \ref{coroloro},

\begin{quote} {\it Let $(S,g,A)$ be an initial data set for the Einstein's Equation,
such that the Riemannian manifold $S$ is simply-connected and
compact, and the tensor field $A$ is negative (resp. positive) definite. Then, there
exist two positive constants $\sigma_i$, $i=1,2$, such that there is no compact spacelike submanifold with future (resp. past) causal mean curvature, on $\phi_{(-\sigma_1,\sigma_2)} S$,
included the extremal case with $\vec{H}=\vec{0}$.

Moreover, that constants are unique and only depend on the initial
value set when their maximality are assumed on the maximal globally
hyperbolic spacetime obtained from the (forwards and backwards)
Cauchy development.}

\end{quote}

Or for example the Theorem \ref{data},

\begin{quote} {\it A spacetime whose timelike sectional curvatures are non-negative
 and which admits an (intrinsic) initial data
set with negative (resp. positive) semi-definite shape operator $A$, does not admit a
submanifold in the future (resp. in the past) of the initial data
set, with future (resp. past)   strictly causal mean curvature vector field.}
\end{quote}

Section \ref{s6} is devoted to the class of standard statics
spacetimes.

Finally, in Sections \ref{s7} and \ref{s8}, some previous results in Section \ref{s6}  are
applied to study some problems of Geometrical Analysis.
Specifically, several Calabi-Bernstein and Dirichlet problems on a
Riemannian manifold are solved.

\vspace{2mm}

\section{Preliminaries}

We denote by $(M,\overline{g})$ an arbitrary $m$-dimensional \emph{
spacetime}, namely, a connected $m$-dimensional oriented and
time-oriented Lorentzian manifold provided with the Lorentzian metric
tensor $\overline{g}$ (see \cite{S-W}).

The notion of symmetry is essential in Physics. In General
Relativity, an infinitesimal symmetry is usually based on the
assumption of the existence of an one-parameter group of
transformations generated by a Killing or, more generally, conformal
Killing (conformal in sort) vector field. In fact, a usual
simplification for the search of exact solutions to the Einstein
equation is to assume the existence, a priori, of such an
infinitesimal symmetry (\cite{Da}, \cite{Er}).

Recall that a vector field $K$ on a spacetime $(M,\overline{g})$ if
called \emph{ conformal} if the Lie derivative of the metric tensor
$\overline{g}$, with respect to $K$ satisfies
$\mathcal{L}_{_K}{\overline{g}}=2\rho{\overline{g}}$, where $\rho$
is a (smooth) function called conformal factor. When $\rho=0$, $K$
is said to be a \emph{ Killing} vector field. If $\rho$ is another
constant, the vector field is called homothetic (see \cite{Da}).

Although different causal characters for the infinitesimal symmetry
may be assumed, the timelike choice is natural, since the integral
curves of a timelike infinitesimal symmetry provide a privileged
class of observers or test particles in the spacetime. Moreover,
this causal choice is supported by very well-known examples of exact
solutions (see, for instance \cite{Da}, \cite{Er}).

\vspace{2mm}

A spacetime $(M,\overline{g})$ admitting a timelike Killing vector
field is called \emph{stationary}. As is well-known, if a spacetime
has a timelike conformal vector field, then it is globally conformal
to a stationary spacetime. For this reason, such spacetimes are
called \emph{conformally stationary} (CS) spacetimes (see,
\cite{ARS}). In general, the orthogonal distribution defined by a
timelike Killing vector field $K$ in a spacetime is not necessarily
integrable. In the special case that this distribution is
integrable, the spacetime is called \emph{static}. A standard model
for a static spacetime is given by a warped product $P\times_h I$,
where $(P,g)$ is a Riemannian manifold, $I$ an open interval and
$h\in C^{\infty}(P)$ a positive smooth function on $P$, endowed with
the metric
\begin{equation}\label{metricagbarra}
\overline{g}=\pi_{_P}^*(g)-h(\pi_{_P})\pi_{_I}^*(dt^2).
\end{equation}
\noindent Here, $\pi_{_P}$  and $\pi_{_I}$ denote the canonical
projections onto the factors $P$ and $I$, respectively (see,
\cite{ON}).

Another important family of spacetimes with an infinitesimal
symmetry are the so called \emph{ pp-wave spacetimes}, i.e.
spacetimes which admit a  parallel global lightlike vector field
\cite[p. 383]{book}.  Moreover, when the spacetime is Ricci-flat
(vacuum solution) it is called \emph{gravitational plane wave}.

\vspace{2mm}

Recall that an isometric immersion $x:S^n \rightarrow M^m$, $2\leq
n<m$, in a spacetime $(M,\overline{g})$  is \emph{spacelike} if the
induced metric via $x$ is Riemannian. In this situation, $S$ is said
to be an (immersed) spacelike submanifold of $M$ (see, \cite[Def.
1.27]{W}). Although from a physical point of view the concept of
submanifold in a spacetime usually corresponds to the case of
embedding submanifold (via the inclusion), the geometric results
that we obtain in this work can be formulated for the more general
framework of immersed submanifolds.
 All the submanifold considered in this work are supposed connected.

The extrinsic geometry of a submanifold $S$ in a spacetime $M$ is
encoded by its second fundamental form $\mathrm{II}:\mathfrak{X}(S)
\times \mathfrak{X}(S) \rightarrow \mathfrak{X}^\bot(S)$, given by
$$
\mathrm{II}(V,W):=\left(\overline{\nabla}_V W\right)^\bot \, ,
$$
where $\overline{\nabla}$ denotes  the Levi-Civita connection of the
metric $\overline{g}$. The \emph{mean curvature vector field} can be
defined according to several conventions. As it is usual in General Relativity, we define it as minus the
metric contraction (without dividing by the dimension of the
submanifold) of the second fundamental form, i.e.,
$$\vec{H}=-\sum_{i=1}^n  \mathrm{II}(E_i,E_i),$$
where $\{E_1,...,E_n\}$ denotes a local orthonormal frame on $S$. Note that this choice of the negative sign is opposite to the one usually taken in Differential Geometry.

\vspace{2mm}

When the immersed submanifold  $x:S \rightarrow M$ is a spacelike
hypersurface (i.e. $m=n+1$) and the induced metric is Riemannian,
the time-orientation of $M$ allows to take a global unitary timelike vector
field $N$ on $S$ pointing to the future.

Let us represent by $g$ the induced metric on the spacelike
hypersurface $S$ and  by $\nabla$ its induced Levi-Civita
connection. The Gauss and Weingarten formulae of $S$ are
respectively
\[
\overline{\nabla}_V W = \nabla_V W - g(AV,W )N \, ,
\]
\[ A V = - \overline{\nabla}_V N \, ,
\] for all
$V, W \in \mathfrak{X}(S)$, where $A$ is the shape operator
associated with $N$. Then, the \textsl{mean curvature function}
associated with $N$ is given, according to our convention, by $H:= 
\mathrm{trace}(A)$.

The mean curvature function is identically  zero if and only if the
spacelike hypersurface is, locally, a critical point of the
$n$-dimensional area functional for compactly supported normal
variations. A spacelike hypersurface with $H=0$ is called a
\textsl{maximal} hypersurface.

As is well-known, spacelike hypersurfaces in general, and of
constant mean curvature in particular, are interesting objects to
study the Einstein equation (see, for instance \cite[Chap.8]{Cho}).
In fact, let $(M,\overline{g})$ be an $(n+1)$-dimensional spacetime
and denote by $\overline{{\rm Ric}}$ and $R(\overline{g})$ its Ricci
tensor and its scalar curvature, respectively. Consider a \emph{
stress-energy tensor field} $T$ on $\overline{M}$, namely a
$2$-covariant symmetric tensor which satisfies some reasonable
conditions from a physical viewpoint (say $T(v,v)\geq 0$ for any
timelike vector $v$, see for example \cite[Section 3.3]{S-W}). It is
said that the spacetime $(M,\overline{g})$ is \emph{an exact
solution to the Einstein equation with zero cosmological constant
and  source $T$},  if the spacetime satisfies
\begin{equation}\label{Einstein}
\overline{{\rm Ric}}-\frac{1}{2}{\rm R}(\overline{g})\overline{g}=T.
\end{equation}

If (\ref{Einstein}) holds, then the following constraint equations are satisfied on each spacelike hypersurface $S$ in $M$
\begin{equation}\label{cond1}
R(g)-{\rm trace}(A^2) + {\rm trace}(A)^2 =\varphi
\end{equation}
\begin{equation}\label{cond2}
{\rm div}(A)-\nabla{\rm trace}(A)= X,
\end{equation}
\noindent where $g$ is the Riemannian metric on $S$ induced by
$\overline{g}$, $R(g)$ its scalar curvature, and $\varphi\in
C^{\infty}(S)$ and $X\in \mathfrak{X}(S)$ depend on the stress
energy tensor $T$. We remark that equations (\ref{cond1}) and
(\ref{cond2}) are respectively obtained from the classical Gauss and
Codazzi equations for the spacelike hypersurface $x:S\rightarrow M$.

Conversely, given $\varphi\in C^{\infty}(S)$ and $X\in
\mathfrak{X}(S)$, they can be seen as differential equations with
unknown $g$ and $A$. Thus, an initial data set for the Cauchy
problem in General Relativity is given by a triple $(S,g,A)$, where
$(S,g)$ is an $n$-dimensional Riemannian manifold and $A:
\mathfrak{X}(S) \rightarrow \mathfrak{X}(S)$ is a $(1,1)$-tensor
field, self-adjoint with respect to $g$, which satisfies the
constraint equations (\ref{cond1}) and (\ref{cond2}).

A solution to the Cauchy problem for the Einstein equation
(\ref{Einstein}) corresponding to the initial data $(S,g,A)$  is a
spacetime $(\overline{M},\overline{g})$, such that $(S,g)$ is an
(embedded) Cauchy spacelike hypersurface and the tensor field $A$
coincides with the shape operator of the embedding. In this setting,
the spacetime $(\overline{M},\overline{g})$ is called a development
of the given initial data set. In \cite{Cho-Ge}, Choquet and Geroch
shown that given an initial data set for the Einstein's equation,
which satisfies the constrain conditions, there exists a
\emph{development of that data set}, which is maximal in the sense
that it is an extension of every other development.

\section{Set up}\label{s3}
Let   $x:S^n \rightarrow M^m$, $n<m$, be an isometric immersed
spacelike submanifold in a spacetime $(M,\overline{g})$.  Let
$\left\{E_i\right\}_{i=1}^n$ and $\left\{N_j\right\}_{j=1}^{m-n}$ be
local orthonormal frames in the tangent vector bundles $TS$ and in
the normal vector bundle $T^\bot S$, respectively. Given a vector
field $X\in\mathfrak{X}(M)$, we define the operator
\begin{eqnarray*}
\mathrm{div}_S X^\bot & := & \sum_i
\overline{g}\left(\overline{\nabla}_{E_i} X^\bot,E_i\right)=
\sum_{i,j} \epsilon_j \, \overline{g}\left(X^\bot,N_j\right) \,
\overline{g}\left(\overline{\nabla}_{E_i} N_j,E_i\right) \\
 & = & \overline{g}\left( X,\vec{H}\right),
\end{eqnarray*}
where $\epsilon_j=\overline{g}(N_j,N_j)$.  Hence,
\begin{equation}\label{div}
\mathrm{div}\left(X^\top\right)= \mathrm{div}_S\left(X\right) -
\overline{g}\left( X,\vec{H}\right),
\end{equation}
where ${\rm div}$ denotes the divergence operator on $(S,g)$, the operator $\mathrm{div}_S$ is defined as
$$\mathrm{div}_S\left(X\right):=\sum_i
\overline{g}\left(\overline{\nabla}_{E_i} X,E_i\right),$$ and
$\vec{H}$ is the mean curvature vector field of $x$. We should point
out that the operator $\mathrm{div}_S$ defined as  (\ref{div}) is
well-defined, i.e., it is independent of the chosen local
orthonormal frame.

In particular, if the submanifold $S$ is compact (without boundary), then making use of the Gauss theorem, we can obtain the following integral formula,
\begin{equation}\label{integral}
\int_S \left\{\mathrm{div}_S\left(X\right) - \overline{g}\left(
X,\vec{H}\right) \right\} \, dV_{g} =0 \,,
\end{equation}
where $dV_g$ is the Riemannian volume element of $(S,g)$.

Here, we must point out several issues. Similar formulas to (\ref{div}) and (\ref{integral}) are obtained in \cite{IB} for the case of a (embedded) parametric spacelike surface in a $4$-dimensional spacetime, Here the authors use a classic coordinate approach.

On the other hand,  in \cite{mars} the authors obtain an intergral formula making use of the variation of the volumen of any subamanifold, so when the variation is determined by a Killing vector field, they can obtain an equality type  (\ref{integral}).

\section{First results}\label{main}

In the first part of this section we will use the integral formula
(\ref{integral}), which will become a powerful tool. We begin with a
general result which allows us to contextualize our framework.

\begin{teor}\label{r1}

 (a) Let $(M,\overline{g})$ be a spacetime, which admits a future causal
vector field $X$ such that the $2$-covariant symmetric tensor
$\mathcal{L}_X \overline{g}$ is positive semi-definite (resp.
negative semi-definite) and positive definite (resp. negative definite) at some point on spacelike vectors. Then there are no
compact spacelike submanifold in $(M,\overline{g})$ with future (resp.
past) causal mean curvature, in particular there are no extremal spacelike submanifold.

\vspace{2mm}

\noindent(b) Let $(M,\overline{g})$ be a spacetime, which admits a future strictly causal
vector field $X$ such that the $2$-covariant symmetric tensor
$\mathcal{L}_X \overline{g}$ is positive semi-definite (resp.
negative semi-definite) on spacelike vectors. Then there are no
compact spacelike submanifolds in $(M,\overline{g})$ with future (resp.
past) causal mean curvature, whose mean curvature vector field
$\vec{H}$ is not zero at least at a point.

\end{teor}

\noindent\emph{Proof.}  (b) Let us suppose that there exists a compact
spacelike submanifold $S$ in $(M,\overline{g})$ with future (resp.
past) causal mean curvature vector field. From the assumption on
the Lie derivative of the metric,   it follows that $\mathrm{div}_S
(X)$ is non-negative (resp. non-positive). Then, since the inner
product of two strictly causal tangent vectors is non-vanishing, the integral
(\ref{integral}) is positive (resp. negative). However, this is not
possible since $S$ is compact without boundary (Gauss theorem). The case (a) is analogous.

\hfill{$\Box$}

\vspace{2mm}

Note that, as particular cases, Theorem \ref{r1} ensurers the
non-existence of weakly, marginally, or nearly future (resp. past)
closed trapped surfaces in such spacetimes. This result must be compared with \cite[Lemma 4.1]{IB}.

%Note that the existence of a future strictly causal vector field is
%equivalent to the existence of a past strictly causal vector.
%
%YO QUITAR\'{I}A ESTE COMENTARIO, NO PEGA MUCHO AQUI
%
%\vspace{2mm}

\vspace{2mm}

We can give a suggestive application of  Theorem \ref{r1} to the
large and relevant class of globally hyperbolic spacetimes. In
\cite{BS}, the authors show that an $(n+1)$-dimensional globally
hyperbolic spacetime $M$ is isometric to a smooth product manifold
\[
\mathbb{R}\times \mathcal{S} \, , \quad \overline{g}= -\beta
d\mathcal{T}^2+\hat{g},
\] where $\mathcal{S}$ is a smooth spacelike
Cauchy hypersurface, $\mathcal{T}: \mathbb{R}\times
\mathcal{S}\rightarrow \mathbb{R}$ the natural projection, $\beta:
\mathbb{R}\times \mathcal{S}\rightarrow (0,\infty)$ a smooth
function, and $\hat{g}$ a $2$-covariant symmetric tensor field on
$\mathbb{R}\times \mathcal{S}$, satisfying:

i) $\nabla{\mathcal{T}}$ is timelike and past-pointing on  $M$.

ii) $\mathcal{S}_{\mathcal{T}}$ is a Cauchy hypersurface  for all
${\mathcal{T}}$, and the restriction $\hat{g}_{_{\mathcal{T}}}$ of
$\hat{g}$ to $\mathcal{S}_{\mathcal{T}}$ is a Riemannian metric.

iii) The radical of $\hat{g}$ at each point $w\in \mathbb{R}\times
\mathcal{S}$ is given by the $Span \{\nabla {\mathcal{T}}\}$ at $w$.

\vspace{2mm}

If we consider the canonical projection $\pi_{_S}$ on $\mathcal{S}$,
which is a diffeomorphism restricted to each level hypersurface
$\mathcal{S}_{\mathcal{T}}$, via $\pi_{_\mathcal{S}}$ we can obtain
an one-parameter family of Riemannian metrics $g_{_\mathcal{T}}$ on
the differentiable manifold $\mathcal{S}$.

In this framework, we can describe a broad family of  spacetimes
(which includes, for instance, the hyperbolic case) by considering a
differentiable manifold $F$, an open interval $I\subset\mathbb{R}$
and an one-parametric family $\{g_t\}_{t\in I}$ of Riemannian
metrics on $F$ (see \cite{ALRUSA} for the details). Then, the
product manifold $M=I\times F$ can be endowed with the Lorentzian
metric given at each point $(t,p)\in I\times F$ by
\[
\overline{g}=\beta\pi_{_\mathbb{R}}^*(-dt^2)+\pi_{_F}^*(g_t) \quad
(\overline{g}=-\beta dt^2+g_t \ \ {\rm in} \ \ {\rm short}),
\]
where $\pi_{_\mathbb{R}}$ and $\pi_{_F}$ denote the canonical
projections onto $I$ and $F$, respectively, and $\beta \in
C^\infty(I \times F)$ is a positive function.  Recall that a
spacetime isometric to a spacetime $(M, \overline{g}=-\beta dt^2+g_t
)$ is called \emph{orthogonal-splitted}. Note that  a spacetime
orthogonal-splitted is stably causal.

It is natural to assume certain homogeneity in the expansive or
contractive behavior of an orthogonal-splitted spacetime. Thus, for
each $v\in T_q F$, $q \in F$, denote by $\tilde{v}$ the lift of $v$
on the integral curve $\alpha_q (s)=(s,q) \in I\times F$. The
spacetime is said to be  \emph{non-contracting in all directions} if
$\partial_t \beta \leq 0$ and $\partial_t
\overline{g}(\tilde{v},\tilde{v})\geq 0$, for all $q\in F$ and $v\in
T_qF$. Analogously, the spacetime is called \emph{non-expanding in
all directions} if the previous inequalities are reversed. Observe
that the spacetime is non-contracting (resp. non-expanding) if and
only if ${\cal L}_{\partial_t}\overline{g}$ is positive
semi-definite (resp. negative semi-definite), where ${\cal L}$
denotes the Lie derivative.

If we consider the observer field
$U=\frac{1}{\sqrt{\beta}}\partial_t$, the proper time $\tau$ of the
observers in $U$ is given by $d\tau=\sqrt{\beta}\,dt$. As a
consequence, the assumption $\partial_t\beta\leq 0$ (resp.
$\partial_t\beta\geq 0$)  means that the rate of change
(acceleration) of the proper time with respect to the time function
$t$ is non-increasing (resp. non-decreasing). Observe also that the
Lorentzian length $\mid\partial_t\mid=-\sqrt{\beta}$ is
non-decreasing (resp. non-increasing) along the integral curves of
$\partial_t$. All  these physical interpretations can be considered
locally.

On the other hand, the spacelike assumption $\partial_t
\overline{g}(\tilde{v},\tilde{v})\geq 0$ (resp. $\partial_t
\overline{g}(\tilde{v},\tilde{v})\leq 0$) guarantees that an
observer in $U$ measures non-contraction (resp. non-expansion) in
all directions of its physical space.

With all of this, as an application of Theorem \ref{r1} to the case
of orthogonal-splitted spacetimes, we have:

\begin{coro}\label{split} Let $(M, \overline{g})$ be an orthogonal-splitted spacetime.
If the spacetime is non-contracting (resp. non-expanding), then
there is no compact spacelike submanifold in $(M,\overline{g})$ with
future (resp. past) causal mean curvature, whose mean curvature
vector field $\vec{H}$ is not zero at least at a point.
\end{coro}

Obviously,  the both Theorem \ref{r1} and Corollary \ref{split} hold
in those regions of the spacetime where the hypotheses on the Lie
derivative of the metric tensor are satisfied.

Furthermore, direct applications to the family of trapped surfaces
can be given from  Corollary \ref{split}.

Note that the causal vector field in Theorem \ref{r1} is not
necessarily associated to an infinitesimal symmetry. When the vector
field is Killing, then the assumptions on the time-orientations can
be avoided. Hence, as a corollary of the Theorem \ref{r1} we can enunciate the following
result, which must be compared with the Theorem 1 in \cite{mars}.

\begin{coro}\label{cor1}
In a spacetime which admits a strictly causal Killing vector field,
there exists no compact spacelike submanifold with future or past
causal mean curvature vector field $\vec{H}$ except, perhaps, with
$\vec{H}\equiv \vec{0}$.

\noindent In particular, in a stationary or pp-wave spacetime, there
exists no compact spacelike submanifold with time-oriented causal
mean curvature vector field except, perhaps, with $\vec{H}\equiv
\vec{0}$.
\end{coro}

Corollary \ref{cor1} has interesting applications to the Cauchy
problem in General Relativity, as becomes clear with the following
result:

\begin{teor}\label{cauchy}
Let $(S,g,A)$ be an initial data set for the Cauchy problem of the
Einstein's equation. Assume that the Riemannian manifold  $S$ admits
a non-minimal compact submanifold $P$, whose mean curvature vector
field $\vec{h}$ in $S$ satisfies
\begin{equation}\label{eqH}
\left\|\vec{h}\right\|< \left|\mathrm{trace}_{|_P} A\right| \, .
\end{equation}
Then, a development spacetime $(M,\tilde{g})$ of this data set
cannot be  a stationary spacetime. In fact, such development cannot
be a spacetime with a strictly causal Killing vector field.
\end{teor}

\noindent\emph{Proof.}

Let us denote by $\nabla^S$ and $\nabla^P$ the induced connections
on $S^n$ and $P^k$ respectively. Let $I\!I_{_S}$, $I\!I_{_P}$ and
$I\!I_{_P}^S$ be the second fundamental forms of $S$ in $M$, of $P$
in ${M}$ and of $P$ in $S$ respectively.

Let us take $\{E_1,...,E_k\}$ a local orthonormal frame on $P$, and
let us extend it to a local orthonormal frame $\{E_1,...E_k,
U_{k+1},...,U_n\}$ on $S$. The mean curvature vector field $\vec{H}$
of $P$ in the spacetime $M$ is given by
\begin{eqnarray*}
-\vec{H}&=&\sum_{i=1}^k
I\!I_{_P}(E_i,E_i)=\sum_{i=1}^k(\overline{\nabla}_{E_i}{E_i})^\perp
\\  & = & \sum_{j=k+1}^n\sum_{i=1}^k
g(\overline{\nabla}_{E_i}{E_i},U_j)U_j-\sum_{i=1}^kg(\overline{\nabla}_{E_i}{E_i},N)N,
\end{eqnarray*}
where $N$ denotes the future unitary normal vector field to $S$ in
the development $M$.

On the other hand, if ${\perp}_{_P}^S$ stands for the orthogonality
to $P$ in $S$,
\begin{eqnarray*}
-\vec{h} & = & \sum_{i=1}^k I\!I_{_P}^S(E_i,E_i)=\sum_{i=1}^k({\nabla}_{E_i}^S E_i)^{{\perp}_{_P}^S} \\
& = &
\sum_{i=1}^k((\overline{\nabla}_{E_i}{E_i})-I\!I_{_S}(E_i,E_i))^{{\perp}_{_P}^S}=\sum_{i=1}^k(\overline{\nabla}_{E_i}{E_i})^{{\perp}_{_P}^S}.
\end{eqnarray*}

Now, taking into account that
\begin{eqnarray*}
\sum_{i=1}^k(\overline{\nabla}_{E_i}{E_i})^\perp=I\!I_{_S}(E_i,E_i)=\sum_{i=1}^k\tilde{g}(\overline{\nabla}_{E_i}{E_i},N)N \\
=\sum_{i=1}^kg(AE_i,E_i)N= \left( \mathrm{trace}_{|_P} A\right),
\end{eqnarray*}
we get that the mean curvature vector field $\vec{H}$ of $P$ in $M$
is related to the mean curvature vector field $\vec{h}$ of $P$ on
$S$ as
\begin{equation}\label{minimal}
\vec{H} = \vec{h} + \left( \mathrm{trace}_{|_P} A\right) \, N \, .
\end{equation}

Observe that $ \left( \mathrm{trace}_{|_P}
A\right)=-\tilde{g}(\vec{H},\vec{N})$ and so, from (\ref{eqH}), we
have that $\vec{H}$ is strictly causal (otherwise, $P$ is minimal in
$S$). Finally, the result follows from Corollary \ref{cor1}.
\hfill{$\Box$}

\begin{rema} {\rm Note that
$\mathrm{trace}_{|_P} A$ is well-defined, i.e., it is independent of
the local frame chosen, as becomes clear from the equality $ \left(
\mathrm{trace}_{|_P} A\right)=-\tilde{g}(\vec{H},\vec{N}).$}
\end{rema}

A more general result than Corollary \ref{cor1} can be given in the
case of existence of a causal conformal vector field. In this
setting, we need ask the conformal factor to be signed:

\begin{coro}\label{con1}
If a spacetime $(M, \overline{g})$ admits a future strictly causal
conformal vector field $K$ whose conformal function $\rho$ is
non-negative (resp. non-positive), then it admits no compact
spacelike submanifold with future (resp. past) time-oriented causal
mean curvature except, perhaps, with  $\vec{H}\equiv 0$.
\end{coro}

\noindent\emph{Proof.} It is enough to observe that ${\rm
div}_S(K)=n\rho$ and to reason as in Theorem \ref{r1}.
\hfill{$\Box$}

If, moreover, we assume that the conformal factor cannot be
identically zero, the conclusion is stronger:

\begin{coro}\label{con2} If a spacetime $(M,\overline{g})$ admits a future strictly causal conformal vector field $K$, whose
conformal factor $\rho$ is a non-zero, non-negative (resp. non-positive) function, then it
admits no compact spacelike submanifold with future (resp. past)
causal mean curvature, including the extremal case with $\vec{H}=\vec{0}$.
\end{coro}

Note that Corollary \ref{con2} applies naturally to spacetimes admitting an homothetic sytrictly causal conformal vector field. Moreover,
both Corollaries \ref{con1} and \ref{con2} admit interesting
applications to the case of (weakly, marginally,.., etc) trapped
surfaces. They also hold for the important family of conformally
stationary spacetimes. So, for instance, we can enunciate (compare with consequences of Lemma 4.1 in \cite{IB}),

\begin{coro} Let $(M,\overline{g})$ be a conformally stationary spacetime, such that the conformal factor of its future timelike conformal vector field is a non-negative and non vanishing identically
(resp. non-positive and non vanishing identically) function. Then,
$M$ does not admit compact future (resp. past) trapped,  nearly
trapped, weakly trapped, marginally trapped and extremal surfaces.

\end{coro}

\section{On the intrinsic Cauchy problem in General Relativity}\label{s5}

%Recall that an initial value set for the
%Cauchy problem in General Relativity is given by a triple $(S,g,A)$, where $(S,g)$ is a Riemannian manifold and $A$ is
%a symmetric operator. A solution for the equations is a spacetime $(\overline{M},\overline{g})$ such that $(S,g)$ is a Cauchy spacelike %hypersurface and $A$ coincides with the shape operator of $S$ in $\overline{M}$.
%We are able to obtain the Cauchy-development of the initial value set.
%The geodesic flow on $S$ will be denoted by $\phi_t$, and we will assume maximality on its definition. Hence, it is easy to see that,
%denoting by $N$ the normal vector field of $S$, we are able to define a new future causal vector field $\overline{N}$ which
%coincides with $N$ at $S$ and such that $\overline{N}|_{\phi_t (p)}=d\phi_t(p) N$, that is, it is a geodesic
%vector field, $\overline{\nabla}_{\overline{N}} \overline{N}=0$. It is easy to see that this vector field is well-defined and unique.
%Moreover, $d\overline{N}^\flat=0$, where $\overline{N}^\flat$ is the one-form associated to $\overline{N}$. Therefore,

In this section, we  study how the existence of a certain compact
(without boundary) spacelike hypersurface in a spacetime
$(M,\overline{g})$ can determine the non-existence of compact
spacelike submanifolds with  strictly causal mean curvature around
$S$. We apply this fact by establishing suitable conditions to a set
of initial data for the Cauchy problem in General Relativity, so
obtaining interesting consequences.

Let $S$ be a spacelike hypersurface immersed in the spacetime $M$.
Assume that $S$ has a global timelike future normal vector field
$N$. We can define a natural extension $\overline{N}$ of ${N}$ on an
open tubular neighborhood of $S$ in $M$. Indeed, given an arbitrary
point $p\in S$ we consider the unique geodesic
$\phi_p(t)=exp_p(tN_p)$ with velocity $N_p$ at $p$. Therefore, the
map $\phi_p(t)$ define the flow of the extension $\overline{N}.$ Due
to the construction itself, it is clear that the $1$-form
$\overline{N}^\flat$ metrically equivalent to $\overline{N}$
satisfies that $d\overline{N}=0$. As a direct consequence we have
that
$$
 \overline{g}\left(\overline{\nabla}_X \overline{N},Y\right) =\frac{1}{2} \left(\mathcal{L}_{\overline{N}} \overline{g}\right)(X,Y)
$$ for $X,Y\in \mathfrak{X}(M)$. Hence, the $(1,1)$-tensor field $\overline{A}: \mathfrak{X}(M) \rightarrow \mathfrak{X}(M)$ defined by
$\overline{g}\left(\overline{A}X,Y\right)=-\frac{1}{2} \left(\mathcal{L}_{\overline{N}} \overline{g}\right)(X,Y)$ extends in a canonical way
the shape operator of $S$.

We denote by $\phi_{(-\epsilon,\epsilon)} S$ the open subset in the
spacetime $M$ given by the points
 $q$ of $M$ such that $q=\phi_p(r)$ for $p\in S$ and $r\in
 (-\epsilon,\epsilon)$. Physically, it is
the portion of the spacetime which is Cauchy-development forwards
and backwards from $S$ up to a quantity of $\epsilon$. Analogously
we can define $\phi_{(-\epsilon,\epsilon)} \Omega$, being $\Omega
\subset S$ a domain in $S$.

\begin{teor}
Let $(S,g,A)$ be an initial data set for the Einstein's Equation
such that the tensor field $A$ is negative (resp. positive) definite. Then, for any
compact domain $\Omega \subset S$, there exists $\epsilon>0$ such
that there is no compact spacelike submanifold with future (resp. past) causal mean curvature, on $\phi_{(-\epsilon,\epsilon)} \Omega$,
included the extremal case with $\vec{H}=\vec{0}$.

\end{teor}

\vspace{2mm}

\noindent\emph{Proof.} Consider a geodesic $\gamma_p :
(-\delta,\delta)\rightarrow M$, with $\gamma_p(0)=p \in \Omega
\subset S$ and $\gamma_p'(0)=N_p$. On $\gamma_p$ we consider the
function $\eta:=\mathrm{det}(\overline{A}\circ \gamma_p)$ 
 Then, there exists a positive number $\delta^+\leq \delta$,
 such that the sign of $\eta$ is constant on $(-\delta^+,\delta^+)$ and this implies
 that $\overline{A}$ is
negative definite on that interval.
Otherwise, it would have a zero eigenvalue and so its determinant would be null, which is a contradiction.

Now, since $\Omega$ is compact, we can find another positive
constant $\epsilon$ such that for any geodesic $\gamma_q:(-\epsilon,
\epsilon)\rightarrow M$, $q \in \Omega$, the same holds.
Equivalently, $\overline{A}$ is negative definite on
$\phi_{(-\epsilon,\epsilon)} \Omega$.

The proof finishes using Theorem \ref{r1} for the timelike vector
field $\overline{N}$, whose divergence is positive.
The other case proof is analogous.
%Assume a time-orientation and time-orient the different geometrical objects.
%For any $\gamma_p:(-\delta,\delta)\rightarrow \overline{M}$, $p\in \Omega$, $\mathrm{det}(A\circ \gamma_p)$ is a function. By assumptions, it %is positive at $0$. Since $\Omega$ is compact, there exists $\epsilon$ such that $\overline{N} \in \mathfrak{X}(\overline{M})$ satisfies
%$\mathcal{L}_{\overline{N}} \overline{g}$ is definite on $\phi_{(-\epsilon,\epsilon)} \Omega$.
%Then, apply Theorem \ref{r1}.
%To end the proof, take a different time-orientation and apply again the same reasoning.
\hfill{$\Box$}

\vspace{2mm}

The assumptions of compactness and $1$-connection can be used to
state the following interesting consequence,

\begin{coro}\label{coroloro}
Let $(S,g,A)$ be an initial data set for the Einstein's Equation,
such that the Riemannian manifold $S$ is simply-connected and
compact, and the tensor field $A$ is negative (resp. positive) definite. Then, there
exist two positive constants $\sigma_i$, $i=1,2$, such that there is no compact spacelike submanifold with future (resp. past) causal mean curvature, on $\phi_{(-\sigma_1,\sigma_2)} S$,
included the extremal case with $\vec{H}=\vec{0}$.

Moreover, that constants are unique and only depend on the initial
value set when their maximality are assumed on the maximal globally
hyperbolic spacetime obtained from the (forwards and backwards)
Cauchy development.
\end{coro}

\begin{rema}
{\rm Corollary \ref{coroloro} admits the following nice topological
interpretation: any compact spacelike submanifold with future (resp. past) causal mean
curvature is far from any simply-connected compact Cauchy
hypersurface with negative (resp. positive) definite shape operator.}
\end{rema}

We can weaken the assumption on the operator $A$, assuming an
initial data set $(S,g,A)$ with $A$ negative semi-definite. % Let us
%look for geometrical conditions under which the thesis in  holds,
%for its extension given by the Lie derivative
%$\mathcal{L}_{\overline{N}}\overline{g}$.
In fact, let $\gamma$ be an arbitrary integral curve of the geodesic
timelike  vector field $\overline{N}$ and consider a spacelike
vector field $X$ on $\gamma$, such that $X$ commute with
$\overline{N}$. Then
\begin{eqnarray*}
-\overline{g}\left(\overline{R}(\overline{N},X)\overline{N},X\right)&=& \overline{g}\left(\overline{\nabla}_{\overline{N}} \overline{\nabla}_X \overline{N},X\right)\\
&=& \frac{1}{2} \gamma'\left( \left(\mathcal{L}_{\overline{N}} \overline{g}\right)(X,X) \right) - \overline{g}\left(\overline{\nabla}_X \overline{N},\overline{\nabla}_X \overline{N} \right) \, .
\end{eqnarray*}
Since the last addend is always non-positive,
%since it is  minus the
%norm of the spacelike vector $\overline{\nabla}_X \overline{N}$.
the previous equality means that $\gamma'\left(
\left(\mathcal{L}_{\overline{N}} \overline{g}\right)(X,X) \right)$
is non-negative if so is the sectional curvature of timelike planes
in the spacetime. In particular, since
$\left(\mathcal{L}_{\overline{N}} \overline{g}\right)(X,X)$ is
non-negative at $S$ (for any tangent vector $X$), then the same
holds in the future of $S$. Thus, if the sectional curvature of
timelike planes is non-negative, then
$\left(\mathcal{L}_{\overline{N}} \overline{g}\right)$ is positive
semi-definite positive in the future of $S$.

Therefore, we can state

\begin{teor}\label{data}
A spacetime whose timelike sectional curvatures are non-negative
 and which admits an (intrinsic) initial data
set with negative (resp. positive) semi-definite shape operator $A$, does not admit a
submanifold in the future (resp. in the past) of the initial data
set, with future (resp. past)   strictly causal mean curvature vector field.
\end{teor}

\begin{rema}
Again, the assumption on the strictly causality of the mean curvature of the
spacelike submanifold can be relaxed. In fact, it suffices to be
time-oriented causal, but not identically zero.
\end{rema}

The techniques used in this section can be applied to obtain
non-existence results of trapped surfaces.

\section{Standard static spacetimes}\label{s6}

In \cite{javasanchez}, under natural causality assumptions, it is
proved that a spacetime $(M,\overline{g})$ which admits a complete
timelike Killing vector field must be standard stationary. That is,
$M$ splits as a topological product of an open interval of the real
line $(\mathbb{R},-dt^2)$ and a Riemannian manifold
$(M_{_0},g_{_0})$ endowed with the metric
$$
\overline{g} = \pi^*_{_0} \, g_{_0} + \pi_{_0}^* \, w \otimes
\pi_{_I}^* dt +  \pi_{_I}^* \, dt \otimes  \pi_{_0}^* \, w -
h(\pi_{_0})^2 \, \pi_{_I}^* dt^2 \, .
$$
Here, $\pi_{_I}$ and $\pi_{_0}$ denote the canonical projections
onto the factors $I$ and $M_{_0}$, respectively, $h$ is a smooth
function on $M_{_0}$ and $w$ is a differential one-form on $M_{_0}$.
If we put $t:=\pi_{_I}$, then the coordinate vector field
$\partial_t:=\frac{\partial}{\partial t}$ is a timelike Killing
vector field and $\overline{\nabla} t= -\frac{1}{h^2} \,
\partial_t$.

A standard stationary spacetime such that the one-form $w$  vanishes
identically is called \textit{standard static}. In this case the
timelike Killing vector field is also irrotational and, as a
consequence, its orthogonal distribution is integrable (see\cite[Ch.
7]{ON}). Following \cite [Ch. 7]{ON}, indeed it is a warped product
that we will denote by $M_{_0} \times_h I$ whose metric can be
written (in short) as
$$\overline{g}=-h^2dt^2+g_{_0}.$$

A relevant family of spacelike hypersurfaces in the standard static
spacetime $M_{_0} \times_h I$ is given by the so called spacelike
slices, i.e., the level hypersurfaces for the time function $t$.
Note that each spacelike slice is a totally geodesic spacelike
hypersurface.

Given a spacelike submanifold $S$ in the spacetime static $M_{_0}
\times_h I$, we will denote the restriction of the temporal function
$t$ to $S$ by $\tau:= \pi_{_I} \circ x$. Then
$$
\nabla \tau = -\frac{1}{h^2} \, \partial_t^\top \, ,
$$ where $\partial_t^\top$ denotes the tangential component
of $\partial_t$ along $S$. Hence, it is not difficult to see that
the Laplacian of $\tau$ on $S$ is given by
\begin{equation}\label{lap1}
\Delta \tau = \frac{2}{h^3} \, \partial_t^\top (h)  + \frac{n}{h^2}
\, \overline{g}\left(\vec{H},\partial_t\right) \, .
\end{equation}

Assume now that $S$ has dimension at least three and consider the
metric
\begin{equation}\label{conformalmetric}
\widetilde{g}=h^{4/(n-2)} g \, ,
\end{equation}
which is conformal to $g$. From (\ref{lap1}) and taking into account
the relation of the Laplacians for conformal metrics (see
\cite{Bes}), it follows that the Laplacian of $\tau$ on $(S,
\widetilde{g})$ is
\begin{equation} \label{conforme}
\widetilde{\Delta} \tau =  n \, h^{-2n/(n-2)} \,
\overline{g}\left(\vec{H},\partial_t\right)  \, .
\end{equation}

\begin{rema}  \label{nota}
{\rm Note that, in the $2$-dimensional case, this conformal metric
cannot be considered. However, next we show a simple procedure which
allows to extend the techniques that we will use, via this conformal
change,   to the $2$-dimensional case.

In fact, given a standard stationary spacetime $(I\times
M_{_0},\overline{g})$, we can consider the new spacetime given by
$(I\times M_{_0} \times \mathbb{S}^2,
\overline{g}+g_{_{\mathbb{S}^2}})$, where $g_{_{\mathbb{S}^2}}$
denotes the usual metric of the Riemannian sphere.  Now, given a
spacelike surface $x:S \rightarrow I\times M_{_0}$, we can define a
new spacelike submanifold in the extended spacetime, $x': S\times
\mathbb{S}^2 \rightarrow  I \times M_{_0} \times \mathbb{S}^2$ as
$x'(p, s)=(x(p),s)$, $p\in S$, $s \in \mathbb{S}^2$. Observe that
the mean curvature vector field of $x'$ is the one of $x$ lifted to
$I \times M_{_0} \times \mathbb{S}^2$, and consequently it keeps the
same causal character. In this setting, we are able to apply the
conformal change described above.

Note that the manifold $\mathbb{S}^2$ may be replaced with
$\mathbb{S}^1$. However, we choose $\mathbb{S}^2$ since it is
simply-connected, and so it does not produce any susceptible change
regarding to the original topologies. Note also that the same
technique can be used for curves (in this case, the mean curvature
is nothing but its acceleration). }
\end{rema}

With all of this, we have:

\begin{teor}\label{teo1}
Let $(M,\overline{g})$ be  a standard static spacetime. Then, there
exists no spacelike submanifold $S$ in $M$ whose mean curvature
vector field $\vec{H}$ satisfies
$\overline{g}\left(\partial_t,\vec{H}\right)\geq 0$ (resp. $\leq 0$)
and such that the function $\tau$ attains a local minimum (resp.
maximum) value on $S$.
\end{teor}

\noindent\emph{Proof.} From
$\overline{g}\left(\partial_t,\vec{H}\right)\geq 0$ and using also
(\ref{conforme}), we have that the function $\tau$  is
$\widetilde{g}$-subharmonic on the compact submanifold $S$.
Therefore, the result easily follows from the classical Hopf maximum
principle (see \cite{Kazdan}, for instance). The case
$\overline{g}\left(\partial_t,\vec{H}\right)\leq 0$ can be reasoned
analogously.

\hfill{$\Box$}

\vspace{2mm}

In the case of compact spacelike submanifolds, we can state the
following rigidity result:

\begin{teor} \label{teo0}
In a standard static spacetime $(M,\overline{g})$, the only compact
spacelike submanifolds such that
$\overline{g}\left(\partial_t,\vec{H}\right) \leq 0$ (or $\geq 0$ ),
must be contained in a spacelike slice.
\end{teor}

\noindent\emph{Proof.} It is enough to observe that the function
$\tau$ is $\tilde{g}$-subharmonic or $\tilde{g}$-superharmonic on a
compact Riemannian manifold, and so $\tau$ must be constant.
\hfill{$\Box$}

The case of vanishing mean curvature corresponds to critical points
of the area functional for normal variations with compact support.
Following the usual nomenclature,  a spacelike submanifold with
vanishing mean curvature is said to be \emph{extremal}.
%They may
%exist, and we obtain the following general characterization for
%them, achieved from previous theorem and making use of standard
%computations,

\begin{coro}\label{extremalcoro}
In a standard static spacetime $(M,\overline{g})$, every extremal
compact submanifold must be contained, as a minimal submanifold, in
a spacelike slice.
\end{coro}

\noindent\emph{Proof.} From Theorem \ref{teo0} we have that an
extremal spacelike submanifold must be contained in a spacelike
slice. Now, taking into account that each slice is totally geodesic
and using also (\ref{minimal}), the result follows. \hfill{$\Box$}

\vspace{2mm}

As regards to the particular case of spacelike hypersurfaces in a
standard static spacetime, note that they admit a globally defined
future-oriented unitary  normal vector field $N$ (where we have
time-oriented the spacetime according to $\partial_t$ pointing to
future). Then, the mean curvature vector field can be written as
$\vec{H}=H \, N$, where $H$ is the \emph{mean curvature function}.
Observe that $\overline{g}\left(N,\frac{\partial_t}{h}
\right)=-\cosh \theta$, where $\theta$ is the hyperbolic angle
between the normal vector field and the observers determined by the
Killing vector field. As another consequence of the conformal change
of metric (\ref{conformalmetric}) and (\ref{conforme}), we can state
the following result, which may be compared with \cite[Theorem
4.6]{ARS}.

\begin{teor} \label{coro1}
In a standard static spacetime  $(M,\overline{g})$, any compact spacelike hypersurface $S$ with signed (in
particular, constant or zero) mean curvature function must be a totally geodesic
spacelike slice.
\end{teor}

\noindent\emph{Proof.} Again, it comes from the subharmonic or
superharmonic character of the function  $\tau$ on the compact
hypersurface $S$. \hfill{$\Box$}

%As commented above, in a standard static spacetime, the orthogonal
%distribution to the timelike Killing vector field $K=\partial_t$ is
%integrable, and its integral hypersurfaces are the spacelike slices
%$\{t=t_{_0}\}$, $t_{_0}\in\mathbb{R}$, which are totally geodesic
%and isometric to the Riemannian manifold (base) $(M,g_{_0})$.

\vspace{2mm}

Although with a different approach,  Theorem  \ref{teo0} and Corollaries \ref{extremalcoro} and \ref{coro1} must be compared with main results in \cite{mars} and \cite{S-S}.

Recall that a complete and simply connected Riemannian manifold is
called a \textit{Cartan-Hadamard} Riemannian manifold provided that
it has non-positive sectional curvatures. In \cite[Cor. 36]{RuSa},
the authors showed that a simply-connected compact manifold cannot
be minimally immersed in a complete Cartan-Hadamard manifold. As a
consequence of this result, we have:

\begin{teor}
Let $(M,\overline{g})$ be a standard static spacetime. Assume that
its base $(M,g_0)$  is a Cartan-Hadamard Riemannian manifold. Then,
the spacetime $M$  admits no  simply-connected compact extremal
submanifolds  with codimension greater than one.
\end{teor}

\noindent\emph{Proof.} It follows from Corollary \ref{extremalcoro}.
\hfill{$\Box$}

The following result is a direct consequence.

\begin{coro}
Let $(S,g,A\equiv 0)$ be an initial value set such that $(S,g)$ is a
Cartan-Hadamard Riemannian manifold. If the spacetime developed is
stationary, then it does not admit a simply-connected compact
extremal submanifold with codimension greater than one.
\end{coro}

In another setting, since the simply-connected compact minimal
surfaces of the round sphere $\mathbb{S}^3$ are totally geodesic
\cite{Bl}, we can state the following characterization of the
maximal surfaces in a highlighted family of spacetimes which
includes the Einstein static spacetime.

\begin{coro}
The only simply-connected compact extremal surfaces  of
$\mathbb{S}^3 \times_h I$ are  the totally geodesic surfaces of
$\mathbb{S}^3 \equiv \left\{ t_{_0} \right\} \times \mathbb{S}^3$,
$t_{_0}\in I$.
\end{coro}

\section{Calabi-Bernstein type problems}\label{s7}

In the literature, the research on maximal hypersurfaces is notably
marked by the discovery of new nonlinear elliptic problems. For
instance, the entire functions defining a minimal graph in the
$(n+1)$-dimensional Lorentz-Minkowski spacetime
${\mathbb{L}}^{n+1}$, satisfy a non-linear elliptic second order
PDE, similar to the equation of minimal graphs in the Euclidean
space $\R^{n+1}$. Nevertheless, in contrast to the Euclidean
setting, in the Lorentzian case the only entire solutions to the
maximal hypersurface equation are the affine functions defining
spacelike hyperplanes. This result was previously proved by Calabi
\cite{Calabi} for $n\leq 4$ and later extended for any $n$ in the
seminal paper by S.Y. Cheng and S.T. Yau \cite{Cheng-Yau}. It is so
known as the Calabi-Bernstein theorem.

Another remarkable difference with the Bernstein theorem for minimal
graphs in the Euclidean space ${\mathbb{R}}^{n+1}$ is that the
Bernstein's Theorem holds only for $n\leq 7$, \cite{S-S-Y}.

In this section, we will  a Calabi-Bernstein type result.
Specifically, we will obtain a uniqueness result for a family of
nonlinear elliptic partial differential equations, which arises from
a geometrical problem for maximal graph in a family of spacetimes.

\vspace{2mm}

Let $(M_0 \times_h I, \overline{g})$ be a standard static spacetime,
whose base is given by the Riemannian $n$-dimensional manifold
$(M_0,g_{_0})$. A function $u\in C^\infty (M_0)$ defines an entire
graph $\Sigma_{u}=\left\{ (p,u(p)) \in M_0 \times_h I \right\}$
 which is an embedded hypersurface in $M_0 \times_h I$.
Given a function $f\in C^\infty(M_0)$, let us denote by $\nabla^{0}
f$ its gradient on $M_0$, that is
 $$\nabla^0 f = \overline{\nabla}
f+\overline{g}(f,\frac{1}{\sqrt{h}} \partial_t)
\frac{\partial_t}{\sqrt{h}}\, .$$
Consequently, the graph
$\Sigma_{u}$ is a spacelike hypersurface if and only if

$$h \,  |\nabla^0 u|^2 <1.$$

On the other hand, the normal vector field of this graph with the
same time-orientation than $-\partial_t$ is
$$
N= \frac{1}{\sqrt{\frac{1}{h}-|\nabla^0 u|^2}} \left\{ \nabla^0 u - \frac{1}{h} \partial_t  \right\} \, .
$$
Then, the hyperbolic angle function of $\Sigma_{u}$ is given by
$$
\cosh \theta = \overline{g}(N,\partial_t)= \frac{1}{\sqrt{1- h \,
|\nabla^0 u|^2}}  \, .
$$
In particular, the graph has bounded hyperbolic angle if and only if
$ h \, |\nabla^0 u|^2$ is bounded away from a certain $\alpha<1$,
which always hols in the spacelike compact case.

If we denote by $N^0$ the projection of the normal vector field $N$
on $TM_0$, then the mean curvature function (i.e. the length of the
mean curvature vector field) of $\Sigma_{u}$ can be computed as
\begin{eqnarray*}
\overline{\mathrm{div}} \,  N & =& \mathrm{div}_0 (N^0) -
\overline{g}\left( \overline{\nabla}_{\frac{\partial_t}{\sqrt{h}}}
N^0, \frac{\partial_t}{\sqrt{h}} \right)
+ \overline{\mathrm{div}}\left( \overline{g}(N,\frac{\partial_t}{\sqrt{h}})\frac{\partial_t}{\sqrt{h}} \right) \\
&=& \mathrm{div}_0 \left(  \frac{ \nabla^0 u }{\sqrt{\frac{1}{h}-|\nabla^0 u|^2}}   \right)
+ {g_{_0}}\left(   \frac{ \nabla^0 u }{\sqrt{\frac{1}{h}-|\nabla^0 u|^2}}, \frac{1}{2} \nabla^0 \log h \right) \, , \\
\end{eqnarray*} where $\mathrm{div}_0$ is the divergence operator of $M_0$.

\vspace{1mm}

We are interested in obtaining all the entire solutions $u$ to the equation
$$\overline{\mathrm{div}} \,  N=H,$$
with the constrain $h \,  |\nabla^0 u|^2 <1$, on a compact
Riemannian manifold,  where $H$ is a known  signed smooth function.
Observe that this constrain assures the elliptic character of the
non-linear problem.

From Theorem \ref{teo0}, we get (compare with the same kind of results presented in \cite{rorusa})
\begin{teor}
Let $(M_0,g_{_0})$ be a compact Riemannian manifold, and let
$h:M_0\longrightarrow \mathbb{R}$ be a positive smooth function.
Consider $H\in C^{\infty}(M_0)$ a signed smooth function. Then, the
equation
$$ \mathrm{div}_0 \left(  \frac{ \nabla^0 u
}{\sqrt{\frac{1}{h}-|\nabla^0 u|^2}}   \right) = \frac{ 1 }{2 \,
\sqrt{\frac{1}{h}-|\nabla^0 u|^2}} \, {g_{_0}}\left( \nabla^0 u,
\nabla^0 \log h  \right)=H
$$  $$\hspace*{4mm} h \, \mid Du\mid^2<1 \, , $$
has solution only provided that $H\equiv 0$. Moreover, in this case
the only solutions are the constant functions.
\end{teor}

\section{Dirichlet problems}\label{s8}

In this final section, we focus on the problem of determining a
spacelike hypersurface with boundary $S$ in a standard  static
spacetime $M_0\times_h\mathbb{R}$, which satisfies that its boundary
is included in a spacelike slice. We will also assume that the mean
curvature function $H$ is signed.

\begin{teor}\label{diri1}
Let $M_0 \times_h \mathbb{R}$ be a standard static spacetime and let
$S$ be an orientable compact hypersurface $S$ whose boundary
$\partial S$ is included in a spacelike slice $\{t=t_{_0}\}$,
$t_0\in\mathbb{R}$. If the mean curvature function $H$ of $S$ is
non-negative (resp. non-positive) and $\tau(S)\geq t_0$ (resp.
$\tau(S)\leq t_{_0}$), then the hypersurface $S$ must be contained
in the slice $\{t=t_{_0}\}$.
\end{teor}

\noindent\emph{Proof.} On the hypersurface $S$, endowed with the
conformal metric $\widetilde{g}$ (and dimensionally extended if
necessary), we consider the vector field $(\tau-t_{_0}) \,
\widetilde{\nabla} \tau$. Obviously,  it vanishes on the boundary
$\partial S$.

Now, we can use the Divergence theorem to obtain
$$
0= \int_S \left(  (\tau-t_{0}) \, \widetilde{\Delta} \tau + |\widetilde{\nabla}\tau|^2_{_{\widetilde{g}}}  \right) \,  dS \, ,
$$ where $dS$ denotes the Riemannian volume element of the manifold $(S,\tilde{g})$. From (\ref{conforme}) it follows that the integrating function
is signed under our assumptions, which implies that $\tau$ must be
constant. \hfill{$\Box$}

%We enunciate in form of corollary the particular case of maximal
%hypersurface.

\begin{coro}
Let $M_0 \times_h I$ be a standard static spacetime. Every
orientable compact maximal hypersurface $S$ such that $\partial S
\subset \left\{t=t_{_0}\right\}$ is contained in the spacelike slice
$\left\{t=t_{_0}\right\}$.
\end{coro}

We can enunciate Theorem \ref{diri1} in term of solutions of a
partial differential inequality as follows

\begin{teor}
Let $\Sigma$ be a compact domain in a Riemannian manifold $(M_0, g_{_0})$ and let $h$ be a smooth positive function defined on $\Sigma$. Then, the only  solutions $u$ defined on $\Sigma$ to the problem
$$
\mathrm{div}_0 \left(  \frac{ \nabla^0 u }{\sqrt{\frac{1}{h}-|\nabla^0 u|^2}}   \right)
- \frac{ 1 }{2 \, \sqrt{\frac{1}{h}-|\nabla^0 u|^2}} \, \overline{g}\left( \nabla^0 u, \nabla^0 \log h  \right)  \leq 0
$$  $$\hspace*{4mm} h \, \mid Du\mid^2<1 \, , $$ $$
u|_{\partial \Sigma} = u_{_0} \, \quad \quad (u_{_0}\in \mathbb{R}) \, ,
$$ $$
u \geq u_{_0} \, ,
$$ are the constant functions.
\end{teor}

Analogously, we have the following Dirichlet problem,

\begin{teor}
Let $\Sigma$ be a compact domain in a Riemannian manifold $(M_0, g_{_0})$ and and let $h$ be a smooth positive function defined on $\Sigma$. Then, the only solutions to the Dirichlet problem on $\Sigma$
$$
\mathrm{div}_0 \left(  \frac{ \nabla^0 u }{\sqrt{\frac{1}{h}-|\nabla^0 u|^2}}   \right)
= \frac{ 1 }{2 \, \sqrt{\frac{1}{h}-|\nabla^0 u|^2}} \, \overline{g}\left( \nabla^0 u, \nabla^0 \log h  \right)
$$  $$\hspace*{4mm} h \, \mid Du\mid^2<1 $$ $$
u|_{\partial \Sigma} = u_{_0} \, \quad \quad (u_{_0}\in \mathbb{R}) \, ,
$$ are the constant functions.
\end{teor}

\section*{Acknowledgments}
Juan A. Aledo was partially supported by MINECO/FEDER Grant no.
MTM2016-80313-P. Rafael M. Rubio and Juan J. Salamanca are partially supported by Spanish MINECO and ERDF Project No. MTM2016-78807-C2-1-P.

\end{document}